\documentclass[12pt]{article}
\usepackage{latexsym, amssymb, amsmath, amscd, amsfonts, epsfig, graphicx, colordvi,verbatim,ifpdf,extarrows}
\usepackage{amsfonts, amsmath, amssymb}
\usepackage{amssymb,amsfonts,amsmath,latexsym,epsfig,cite, psfrag,eepic,color}
\usepackage{amscd,graphics}
\usepackage{latexsym, amssymb,  amsmath,amscd, amsfonts, epsfig, graphicx, colordvi,amsthm}

\usepackage{graphicx}
\usepackage{epstopdf}
\usepackage{color}
\usepackage{ifpdf}
\usepackage{fancybox}
\usepackage[font=small,labelfont=bf,labelsep=none]{caption}
\usepackage{float}

\allowdisplaybreaks

\newtheorem{thm}{Theorem}[section]

\newtheorem{defi}[thm]{Definition}
\newtheorem{lem}[thm]{Lemma}
\newtheorem{core}[thm]{Corollary}

\def\pf{\noindent{\it Proof.} }
\def\qed{\nopagebreak\hfill{\rule{4pt}{7pt}}
\medbreak}

\numberwithin{equation}{section}

\def\qed{\nopagebreak\hfill{\rule{4pt}{7pt}}
\medbreak}

\newlength{\boxedparwidth}
  {\begin{center} \begin{tabular}{|@{\hspace{.315in}}c@{\hspace{.15in}}|}
                  \hline \\ \begin{minipage}[t]{\boxedparwidth}
                  \setlength{\parindent}{.25in}}%
  {\end{minipage} \\ \\ \hline \end{tabular} \end{center}}

\begin{document}

\begin{center}
{\Large \bf Separable integer partition classes with restrictions on consecutive parts}
\end{center}

\begin{center}
 {Y.Q. Chen}$^{1}$, {Thomas Y. He}$^{2}$, {X.M. Huang}$^{3}$ and
  {T.T. Zou}$^{4}$ \vskip 2mm

$^{1,2,3,4}$ School of Mathematical Sciences, Sichuan Normal University, Chengdu 610066, P.R. China

   \vskip 2mm

  $^1$chenyuqian11@stu.sicnu.edu.cn, $^2$heyao@sicnu.edu.cn,  $^3$huangxinmin@stu.sicnu.edu.cn, $^4$zoutingting@stu.sicnu.edu.cn
\end{center}

\vskip 6mm   {\noindent \bf Abstract.} Recently, Andrews introduced separable integer partition classes and studied some well-known theorems. In this article, we will consider the types of partitions with restrictions on consecutive parts. We will show that such partitions are separable integer partition classes and then give the generating functions for such partitions.



\section{Introduction}

A partition $\pi$ of a positive integer $n$ is a finite non-increasing sequence of positive integers $\pi=(\pi_1,\pi_2,\ldots,\pi_m)$ such that $\pi_1+\pi_2+\cdots+\pi_m=n$. The empty sequence forms the only partition of zero. The $\pi_i$ are called the parts of $\pi$. The weight of $\pi$ is the sum of parts, denoted $|\pi|$. Assume that $k$ and $d$ are integers such that $k\geq d\geq 1$, we use $\ell_{k,d}(\pi)$ to denote the number of parts congruent to $d$ modulo $k$ in $\pi$.

In \cite{Andrews-2022}, Andrews introduced separable integer partition classes and analyzed some well-known theorems, such as the first G\"ollnitz-Gordon identity, Schur's partition theorem, partitions  with $n$ copies of $n$, and so on.

 \begin{defi}\cite[Page 619]{Andrews-2022}
For a positive integer $k$, a separable integer partition class $\mathcal{P}$ with modulus $k$ is a subset of all the partitions satisfying the following{\rm:}

There is a subset $\mathcal{B}\subset\mathcal{P}$ {\rm(}$\mathcal{B}$ is called the basis of $\mathcal{P}${\rm)} such that for each integer $m\geq 1$, the number of partitions in $\mathcal{B}$ with $m$ parts is finite and every partition in $\mathcal{P}$ with $m$ parts is uniquely of the form
\begin{equation}\label{ordinary-form-1}
(b_1+\pi_1)+(b_2+\pi_2)+\cdots+(b_m+\pi_m),
\end{equation}
where $(b_1,b_2,\ldots,b_m)$ is a partition in $\mathcal{B}$ and $(\pi_1,\pi_2,\ldots,\pi_m)$ is a  non-increasing sequence of nonnegative integers, whose only restriction is that each part is divisible by $k$. Furthermore, all partitions of the form \eqref{ordinary-form-1} are in $\mathcal{P}$.
\end{defi}

Assume that $\mathcal{P}$ is a separable integer partition class with modulus $k$ and  $\mathcal{B}$ is the basis of $\mathcal{P}$. For $m\geq 1$, let $\mathcal{B}(m)$ be the set of partitions in $\mathcal{B}$ with  $m$ parts. Clearly, the generating function for the partitions in $\mathcal{P}$ is
 \begin{equation*}
  \sum_{\pi\in\mathcal{P}}z^{\ell_{k,d}(\pi)}q^{|\pi|}=1+\sum_{m\geq 1}\frac{1}{(q^k;q^k)_m}\sum_{\pi\in\mathcal{B}(m)}z^{\ell_{k,d}(\pi)}q^{|\pi|}.
  \end{equation*}
  
 Here and in the sequel, we assume that $|q|<1$. The $q$-shifted factorial is denoted by
\[(a;q)_n=\left\{\begin{array}{ll}
1,&\text{if }n=0,\\
(1-a)(1-aq)\cdots(1-aq^{n-1}),&\text{if }n\geq 1,
\end{array}\right.\]
For positive integer $k$, the $q$-binomial coefficient, or Gaussian polynomial for integers $A$ and $B$ is given by
\[{A\brack B}_k=\left\{\begin{array}{ll}\frac{(q^k;q^k)_A}{(q^k;q^k)_B(q^k;q^k)_{A-B}},&\text{if }A\geq B\geq 0,\\
0,&\text{otherwise,}
\end{array}
\right.\]
which is the generating function for the partitions such that there are at most $B$ parts and each part is divisible by $k$ and does not exceed $k(A-B)$.

There are several works related to separable integer partition classes, please refer to \cite{Chen-He-Hu-Xie-2024,Chen-He-Huang-Zhang-2025,Chen-He-Tang-Wei-2024,He-Hu-Huang-Xie-2025,He-Huang-Li-Zhang-2025,Passary-2019}.
Especially, with the aid of separable integer partition classes, He, Huang, Li and Zhang studied $\mathcal{P}_{a,b,k}$ in \cite{He-Huang-Li-Zhang-2025} and  He, Hu, Huang and Xie investigated $\mathcal{R}_{a,b,k}$ in \cite{He-Hu-Huang-Xie-2025} for $1\leq a<b\leq k$, where $\mathcal{P}_{a,b,k}$ is the set of partitions whose parts are congruent to $a$ or $b$ modulo $k$ and $\mathcal{R}_{a,b,k}$ is the set of partitions in $\mathcal{P}_{a,b,k}$ such that only parts congruent to $a$ modulo $k$ may be repeated.

One of the main objectives of this article is to study four types of partitions from the point of view of separable integer partition classes. We first define two subsets $\mathcal{P}_{a,b,k,r}$  and $\mathcal{P}'_{a,b,k,r}$ of $\mathcal{P}_{a,b,k}$.
 \begin{defi} Let $k,a,b$ and $r$ be integers such that $k\geq b>a\geq 1$ and $r\geq 1$. We define $\mathcal{P}_{a,b,k,r}$ (resp. $\mathcal{P}'_{a,b,k,r}$) to be the set of partitions such that
 \begin{itemize}
 \item[(1)] the parts are congruent to $a$ or $b$ modulo $k$;
 \item[(2)] there are no $r+1$ consecutive parts which are congruent to $b$ (resp. $a$) modulo $k$.
 \end{itemize}
\end{defi}

For example, there are eleven partitions of $7$ in $\mathcal{P}_{1,2,2,1}$.
\[(7),(6,1),(5,2),(5,1,1),(4,3),(4,1,1,1),(3,3,1),\]
\[(3,2,1,1),(3,1,1,1,1),(2,1,1,1,1,1),(1,1,1,1,1,1,1).\]
There are eight partitions of $7$ in $\mathcal{P}'_{1,2,2,2}$.
\[(7),(6,1),(5,2),(4,3),(4,2,1),(3,2,2),(3,2,1,1),(2,2,2,1).\]

Then, we give the definitions of $\mathcal{R}_{a,b,c,k}$ and $\mathcal{R}_{a,b,c,k,r}$.
 \begin{defi} Assume that $k,a,b,c$ and $r$ are integers such that $k\geq c>b>a\geq 1$ and $r\geq 1$. Let $\mathcal{R}_{a,b,c,k}$ be the set of partitions $\pi=(\pi_1,\pi_2,\ldots,\pi_m)$ such that
 \begin{itemize}
 \item[(1)] the parts are congruent to $a$ or $b$ or $c$ modulo $k$;
 \item[(2)] the parts congruent to $c$ modulo $k$ are distinct;
 \item[(3)] for $1\leq i<m$, if $\pi_i\equiv a\pmod{k}$ then $\pi_{i+1}$ is congruent to $a$ or $c$ modulo ${k}$.
 \end{itemize}
 Let $\mathcal{R}_{a,b,c,k,r}$ be the set of partitions in $\mathcal{R}_{a,b,c,k}$ such that there are no $r$ consecutive parts which are congruent to $b$ modulo $k$.
 \end{defi}

For example, there are nine partitions of $12$ in $\mathcal{R}_{2,3,4,4,2}$.
\[(12),(10,2),(8,4),(8,2,2),(6,6),(6,4,2),(6,2,2,2),(4,2,2,2,2),(2,2,2,2,2,2).\]
Besides, there are four partitions of $12$ belonging to $\mathcal{R}_{2,3,4,4}$ but not in  $\mathcal{R}_{2,3,4,4,2}$.
\[(7,3,2),(4,3,3,2),(3,3,3,3),(3,3,2,2,2).\]

Clearly, if $r=1$, then the set $\mathcal{R}_{a,b,c,k,1}$ reduces to $\mathcal{R}_{a,c,k}$.

We will show that $\mathcal{P}_{a,b,k,r}$, $\mathcal{P}'_{a,b,k,r}$, $\mathcal{R}_{a,b,c,k}$ and $\mathcal{R}_{a,b,c,k,r}$ are separable integer partition classes with modulus $k$ and give the generating functions for the partitions in $\mathcal{P}_{a,b,k,r}$, $\mathcal{P}'_{a,b,k,r}$, $\mathcal{R}_{a,b,c,k}$ and $\mathcal{R}_{a,b,c,k,r}$.
 \begin{thm}\label{thm-p-abkr-0} For $k\geq b>a\geq 1$ and $r\geq 1$,
\begin{equation}\label{gen-p-abkr-0}
 \begin{split}
 &\quad\sum_{\pi\in\mathcal{P}_{a,b,k,r}}\mu^{\ell_{k,a}(\pi)}\nu^{\ell_{k,b}(\pi)}q^{|\pi|}\\
 &=1+\sum_{m\geq 1}\frac{1}{(q^k;q^k)_m}\left(\mu^mq^{ma}+\sum_{s\geq 1}\sum_{h\geq 0}\mu^{m-h-s}\nu^{h+s}q^{(m-h-s)a+(h+s)b+k(s^2-s)}{{m-h-s+1}\brack {s}}_{k}\right.\\
 &\left.\qquad\qquad\qquad\qquad\qquad\qquad\qquad\times\sum_{i\geq 0}(-1)^iq^{rk(i^2-i)/2}{{s}\brack {i}}_{rk}{{h-ri+s-1}\brack {s-1}}_{k}\right),
 \end{split}
 \end{equation}
 and
  \begin{equation}\label{gen-pp-abkr-0}
 \begin{split}
&\quad\sum_{\pi\in\mathcal{P}'_{a,b,k,r}}\mu^{\ell_{k,a}(\pi)}\nu^{\ell_{k,b}(\pi)}q^{|\pi|}\\
&=1+\sum_{m\geq 1}\frac{1}{(q^k;q^k)_m}\left\{\nu^mq^{mb}+\sum_{s\geq 1}\sum_{h\geq 0}\mu^{h+s}\nu^{m-h-s}q^{(h+s)a+(m-h-s)b+k(s-1)^2}\right.\\
&\times\left.\left({{m-h-s}\brack {s-1}}_{k}+q^{k(h+2s-1)}{{m-h-s}\brack {s}}_{k}\right)\sum_{i\geq 0}(-1)^iq^{rk(i^2-i)/2}{{s}\brack {i}}_{rk}{{h-ri+s-1}\brack {s-1}}_{k}\right\}.
 \end{split}
 \end{equation}
\end{thm}

 \begin{thm}\label{thm-r-abkr-0} For $k\geq c>b>a\geq 1$ and $r\geq 1$,
 \begin{align*}
&\quad\sum_{\pi\in\mathcal{R}_{a,b,c,k}}\mu^{\ell_{k,a}(\pi)}\nu^{\ell_{k,b}(\pi)}\omega^{\ell_{k,c}(\pi)}q^{|\pi|}\\
&=1+\sum_{m\geq 1}\frac{1}{(q^k;q^k)_m}\sum_{s,h\geq 0}\mu^{m-h-s}\nu^{h}\omega^sq^{(m-h-s)a+hb+sc+k(s^2-s)/2}{{h+s}\brack {s}}_{k}{{m-h}\brack {s}}_{k},
 \end{align*}
 and
\begin{align*}
&\quad\sum_{\pi\in\mathcal{R}_{a,b,c,k,r}}\mu^{\ell_{k,a}(\pi)}\nu^{\ell_{k,b}(\pi)}\omega^{\ell_{k,c}(\pi)}q^{|\pi|}\\
&=1+\sum_{m\geq 1}\frac{1}{(q^k;q^k)_m}\sum_{s,h\geq 0}\mu^{m-h-s}\nu^{h}\omega^sq^{(m-h-s)a+hb+sc+k(s^2-s)/2}{{m-h}\brack {s}}_{k}\\
&\qquad\qquad\qquad\qquad\qquad\times\sum_{i\geq 0}(-1)^iq^{rk(i^2-i)/2}{{s+1}\brack {i}}_{rk}{{h-ri+s}\brack {s}}_{k}.
\end{align*}
 \end{thm}

In this article, we also investigate the overpartitions with restrictions on consecutive parts.
An overpartition, introduced by Corteel and Lovejoy \cite{Corteel-Lovejoy-2004},  is a partition such that the first occurrence (equivalently, the last occurrence) of a number may be overlined. Let $\pi$ be an overpartition. We use $\ell_o(\pi)$ to denote the number of overlined parts in $\pi$. For a part $\pi_i$ of $\pi$, we set $|\pi_i|=t$  if $\pi_i=t$ or $\overline{t}$. We adopt the following convention: For positive integer $t$ and nonnegative integer $d$, we define $\overline{t}+d=\overline{t+d}$.

In \cite{Chen-He-Tang-Wei-2024}, Chen, He, Tang and Wei extended separable integer partition classes with modulus $1$ introduced by Andrews \cite{Andrews-2022} to overpartitions, called separable overpartition classes, and studied overpartitions and the overpartition analogue of Rogers-Ramanujan identities in view of separable overpartition classes.

\begin{defi}\cite[Definition 4.1]{Chen-He-Tang-Wei-2024}\label{defi-separable}
A separable overpartition class $\mathcal{P}$ is a subset of all the overpartitions satisfying the following{\rm:}

There is a subset $\mathcal{B}\subset\mathcal{P}$ {\rm(}$\mathcal{B}$ is called the basis of $\mathcal{P}${\rm)} such that for each integer $m\geq 1$, the number of overpartitions in $\mathcal{B}$ with $m$ parts is finite and every overpartition in $\mathcal{P}$ with $m$ parts is uniquely of the form
\begin{equation}\label{over-form-1}
(b_1+\pi_1)+(b_2+\pi_2)+\cdots+(b_m+\pi_m),
\end{equation}
where $(b_1,b_2,\ldots,b_m)$ is an overpartition in $\mathcal{B}$ and $(\pi_1,\pi_2,\ldots,\pi_m)$ is a  non-increasing sequence of nonnegative integers. Moreover, all overpartitions of the form \eqref{over-form-1} are in $\mathcal{P}$.
\end{defi}

Assume that $\mathcal{P}$ is a separable overpartition class and  $\mathcal{B}$ is the basis of $\mathcal{P}$. For $m\geq 1$ and $s\geq 0$, let $\mathcal{B}(m,s)$ be the set of overpartitions in $\mathcal{B}$ with $m$ parts and $s$ overlined parts. Then, the generating function for  the overpartitions in $\mathcal{P}$ is
 \begin{equation*}
  \sum_{\pi\in\mathcal{P}}z^{\ell_o(\pi)}q^{|\pi|}=1+\sum_{m\geq 1}\sum_{s\geq 0}\frac{z^s}{(q;q)_m}\sum_{\pi\in\mathcal{B}(m,s)}q^{|\pi|}.
  \end{equation*}

In this article, we will study four types of overpartitions from the point of view of separable overpartition classes.
\begin{thm}\label{main-over}
Let $\overline{\mathcal{F}}$ (resp. $\overline{\mathcal{L}}$) be the set of overpartitions in which the first (resp. last) occurrence of a number may be overlined and no adjacent parts are both overlined. For $r\geq 1$, let $\mathcal{F}_r$ (resp. $\mathcal{L}_r$) be the set of overpartitions in which the first (resp. last) occurrence of a number may be overlined and
there are no $r$ consecutive parts which are non-overlined. Then,
\begin{equation}\label{eqn-over-f-over}
\sum_{\pi\in\overline{\mathcal{F}}}z^{\ell_o(\pi)}q^{|\pi|}=1+\sum_{m\geq 1}\sum_{s\geq 0}\frac{z^sq^{m+s^2-s}}{(q;q)_m}{{m-s+1}\brack{s}}_1,
\end{equation}
\begin{equation}\label{eqn-over-L-over}
\sum_{\pi\in\overline{\mathcal{L}}}z^{\ell_o(\pi)}q^{|\pi|}=1+\sum_{m\geq 1}\sum_{s\geq 0}\frac{z^sq^{m+(s-1)^2}}{(q;q)_m}\left({{m-s}\brack{s-1}}_1+q^{2s-1}{{m-s}\brack{s}}_1\right),
\end{equation}
\begin{equation}\label{eqn-over-f-non-over}
\sum_{\pi\in\mathcal{F}_r}z^{\ell_o(\pi)}q^{|\pi|}=1+\sum_{m\geq 1}\sum_{s\geq 0}\frac{z^sq^{m+(s^2-s)/2}}{(q;q)_m}\sum_{i\geq 0}(-1)^iq^{r(i^2-i)/2}{{s+1}\brack {i}}_{r}{{m-ri}\brack {s}}_{1},
\end{equation}
and
\begin{equation}\label{eqn-over-L-non-over}
\begin{split}
&\quad\sum_{\pi\in\mathcal{L}_r}z^{\ell_o(\pi)}q^{|\pi|}\\
&=1+\sum_{m\geq 1}\frac{q^m}{(q;q)_m}\left\{{{r-m-1}\brack {0}}_{1}+\sum_{s\geq 1}z^sq^{m+(s^2-s)/2}\left(\sum_{i\geq 0}(-1)^iq^{r(i^2-i)/2}{{s}\brack {i}}_{r}{{m-ri-1}\brack {s-1}}_{1}\right.\right.\\
&\qquad\qquad\qquad\qquad\qquad\left.\left.+q^m\sum_{j=1}^{r-1}\sum_{i\geq 0}(-1)^iq^{r(i^2-i)/2-j}{{s}\brack {i}}_{r}{{m-j-ri-1}\brack {s-1}}_{1}\right)\right\}.
\end{split}
\end{equation}
\end{thm}

For example, there are twelve overpartitions of $4$ in $\overline{\mathcal{F}}$.
\[(4),(\overline{4}),(3,1),(\overline{3},1),(3,\overline{1}),(2,2),(\overline{2},2),\]
\[(2,1,1),(\overline{2},1,1),(2,\overline{1},1),(1,1,1,1),(\overline{1},1,1,1).
\]
There are thirteen overpartitions of $4$ in $\overline{\mathcal{L}}$.
\[(4),(\overline{4}),(3,1),(\overline{3},1),(3,\overline{1}),(2,2),({2},\overline{2}),\]
\[(2,1,1),(\overline{2},1,1),(2,1,\overline{1}),(\overline{2},1,\overline{1}),(1,1,1,1),(1,1,1,\overline{1}).\]
There are eight overpartitions of $4$ in $\mathcal{F}_2$.
\[(4),(\overline{4}),(\overline{3},1),(3,\overline{1}),(\overline{3},\overline{1}),(\overline{2},2),(2,\overline{1},1),(\overline{2},\overline{1},1).\]
There are seven overpartitions of $4$ in $\mathcal{L}_2$.
\[(4),(\overline{4}),(\overline{3},1),(3,\overline{1}),(\overline{3},\overline{1}),({2},\overline{2}),(\overline{2},1,\overline{1}).\]

This article is organized as follows. In section 2, we first investigate the function
\[g_{k,r}(h,s)=\sum_{i\geq 0}(-1)^iq^{rk(i^2-i)/2}{{s}\brack {i}}_{rk}{{h-ri+s-1}\brack {s-1}}_{k},\]
where $k\geq 1$, $r\geq 1$, $h\geq 0$ and $s\geq 1$. Then, we will show Theorems \ref{thm-p-abkr-0}, \ref{thm-r-abkr-0} and \ref{main-over} in Section 3, Section 4 and Section 5 respectively.

\section{The function $g_{k,r}(h,s)$}

In this section, we assume that $k,d,r,h$ and $s$ are integers such that $k\geq d\geq 1$, $r\geq 1$, $h\geq 0$ and $s\geq 1$. The objective of this section is to investigate the function $g_{k,r}(h,s)$.
\begin{defi}
Let $\mathcal{G}_{d,k,r}(h,s)$ be the set of partitions such that
\begin{itemize}
\item[(1)] there are $h$ parts;
\item[(2)] all parts are congruent to $d$ modulo $k$ and do not exceed $k(s-1)+d$;
\item[(3)] each part appears at most $r-1$ times.
\end{itemize}
\end{defi}
The following theorem plays an important role in the proofs of Theorems \ref{thm-p-abkr-0}, \ref{thm-r-abkr-0} and \ref{main-over}.

\begin{thm}\label{useful-gen-h-s}
The generating function for the partitions in $\mathcal{G}_{d,k,r}(h,s)$ is
\[\sum_{\pi\in\mathcal{G}_{d,k,r}(h,s)}q^{|\pi|}=q^{hd}g_{k,r}(h,s)=q^{hd}\sum_{i\geq 0}(-1)^iq^{rk(i^2-i)/2}{{s}\brack {i}}_{rk}{{h-ri+s-1}\brack {s-1}}_{k}.\]
\end{thm}

Before proving Theorem \ref{useful-gen-h-s}, we recall two identities due to Cauchy \cite[p. 46]{Cauchy-1893}, see also \cite[Theorem 3.3]{Andrews-1976}).
\begin{equation}\label{Cauchy-1}
(z;q)_s=\sum_{i\geq 0}(-1)^iz^iq^{(i^2-i)/2}{{s}\brack{i}}_1,
\end{equation}
and
\begin{equation}\label{Cauchy-2}
\frac{1}{(z;q)_s}=\sum_{i\geq 0}z^i{{i+s-1}\brack{s-1}}_1.
\end{equation}

Now, we proceed to show Theorem \ref{useful-gen-h-s}.

{\noindent \bf Proof of Theorem \ref{useful-gen-h-s}.} It is clear that
\begin{align}
\sum_{h\geq 0}z^h\sum_{\pi\in\mathcal{G}_{d,k,r}(h,s)}q^{|\pi|}&=\prod_{i=0}^{s-1}\left(1+zq^{ik+d}+z^2q^{2(ik+d)}+\cdots+z^{r-1}q^{(r-1)(ik+d)}\right)\nonumber\\
&=\prod_{i=0}^{s-1}\frac{1-z^rq^{r(ik+d)}}{1-zq^{ik+d}}\nonumber\\
&=\frac{(z^rq^{rd};q^{rk})_s}{(zq^b;q^k)_s}.\label{inter-cauchy}
\end{align}

Letting $q\rightarrow q^{rk}$ and then $z\rightarrow z^{r}q^{rd}$ in \eqref{Cauchy-1}, we have
 \begin{equation}\label{use-Cauchy-1}
(z^{r}q^{rd};q^{rk})_s=\sum_{i\geq 0}(-1)^i(z^{r}q^{rd})^iq^{rk(i^2-i)/2}{{s}\brack{i}}_{rk}.
\end{equation}

Letting $q\rightarrow q^{k}$ and then $z\rightarrow zq^{d}$ in \eqref{Cauchy-2}, we have
\begin{equation}\label{use-Cauchy-2}
\frac{1}{(zq^{d};q^k)_s}=\sum_{i\geq 0}(zq^d)^i{{i+s-1}\brack{s-1}}_{k}.
\end{equation}

Substituting \eqref{use-Cauchy-1} and \eqref{use-Cauchy-2} into \eqref{inter-cauchy}, we get
\begin{align*}
\sum_{h\geq 0}z^h\sum_{\pi\in\mathcal{G}_{d,k,r}(h,s)}q^{|\pi|}&=\sum_{i,j\geq 0}(-1)^iz^{ri+j}q^{rk(i^2-i)/2+(ri+j)d}{{s}\brack{i}}_{rk}{{j+s-1}\brack{s-1}}_{k}\\
&=\sum_{h\geq 0}z^hq^{hd}\sum_{i\geq 0}(-1)^iq^{rk(i^2-i)/2}{{s}\brack{i}}_{rk}{{h-ri+s-1}\brack{s-1}}_{k}.
\end{align*}

By comparing the coefficients of $z^h$ on the extremes of the string of equations above, we complete the proof.   \qed

For a partition $\pi$ in $\mathcal{G}_{d,k,r}(h,s)$, the number of parts in $\pi$ is at most $(r-1)s$. This implies that for $h>(r-1)s$, the set $\mathcal{G}_{d,k,r}(h,s)$ is empty. The following corollary immediately follows from  Theorem \ref{useful-gen-h-s}.
\begin{core}
For $h>(r-1)s$, we have
\[g_{k,r}(h,s)=\sum_{i\geq 0}(-1)^iq^{rk(i^2-i)/2}{{s}\brack {i}}_{rk}{{h-ri+s-1}\brack {s-1}}_{k}=0.\]
\end{core}

It is easy to get the following generating function for the partitions in $\mathcal{G}_{d,k,r}(h,s)$ for $r=2$.
\begin{thm}\label{useful-gen-h-s-d=2}
\[\sum_{\pi\in\mathcal{G}_{d,k,2}(h,s)}q^{|\pi|}=q^{hd+k(h^2-h)/2}{{s}\brack {h}}_{k}.\]
\end{thm}

\pf  Clearly, we have
\begin{equation*}
\sum_{h\geq 0}z^h\sum_{\pi\in\mathcal{G}_{d,k,2}(h,s)}q^{|\pi|}=\prod_{i=0}^{s-1}\left(1+zq^{ik+d}\right)=(-zq^d;q^k)_s=\sum_{h\geq 0}z^hq^{hd+k(h^2-h)/2}{{s}\brack {h}}_{k},
\end{equation*}
where the final equation follows from \eqref{Cauchy-1} with $q\rightarrow q^k$ and then $z\rightarrow -zq^d$. Comparing the coefficients of $z^h$ on the extremes of the string of equations above, we complete the proof.   \qed

Combining Theorems \ref{useful-gen-h-s} and \ref{useful-gen-h-s-d=2}, we can get the following corollary.
\begin{core}
\[\sum_{i\geq 0}(-1)^iq^{k(i^2-i)}{{s}\brack {i}}_{2k}{{h-2i+s-1}\brack {s-1}}_{k}=q^{k(h^2-h)/2}{{s}\brack {h}}_{k}.\]
\end{core}

\section{Proof of Theorem \ref{thm-p-abkr-0}}

In this section, we assume that $k,a,b$ and $r$ are integers such that $k\geq b>a\geq 1$ and $r\geq 1$. In order to prove Theorem \ref{thm-p-abkr-0},  we will show that  $\mathcal{P}_{a,b,k,r}$  and $\mathcal{P}'_{a,b,k,r}$ are separable integer partition classes with modulus $k$ and then give the proofs of \eqref{gen-p-abkr-0} and \eqref{gen-pp-abkr-0} in Section 3.1 and Section 3.2 respectively.

\subsection{Proof of \eqref{gen-p-abkr-0}}

We first show that $\mathcal{P}_{a,b,k,r}$ is a separable integer partition class with modulus $k$. To do this, we are required to find the basis of $\mathcal{P}_{a,b,k,r}$, which involves the following set.
\begin{defi}\label{basis-p-m}
For $m\geq 1$, let $\mathcal{BP}_{a,b,k,r}(m)$ denote the set of partitions $\pi=(\pi_1,\pi_2,\ldots,\pi_m)$ such that
\begin{itemize}
 \item[(1)] for $1\leq i\leq m$, $\pi_i$ is congruent to $a$ or $b$ modulo ${k}$;
 \item[(2)] $\pi_m=a$ or $b$;
 \item[(3)] for $1\leq i<m$, $\pi_i-\pi_{i+1}<k$;
 \item[(4)] for $1\leq i\leq m-r$, if $\pi_{i+1}=\pi_{i+2}=\cdots=\pi_{i+r}\equiv b\pmod{k}$ then $\pi_i\equiv a\pmod{k}$.
 \end{itemize}
\end{defi}

Under the conditions (1)-(3) in Definition \ref{basis-p-m}, we find that for $m\geq 1$, the number of partitions in $\mathcal{BP}_{a,b,k,r}(m)$ does not exceed $2^m$. For example, there are thirteen partitions in $\mathcal{BP}_{a,b,k,2}(4)$.
\[(a,a,a,a),(b,a,a,a),(b,b,a,a),(k+a,b,a,a),(k+a,b,b,a),\]
\[(k+a,k+a,b,a),(k+b,k+a,b,a),(k+a,k+a,b,b),(k+b,k+a,b,b),\]
\[(k+a,k+a,k+a,b),(k+b,k+a,k+a,b),(k+b,k+b,k+a,b),(2k+a,k+b,k+a,b).\]

Set
\[\mathcal{BP}_{a,b,k,r}=\bigcup_{m\geq1}\mathcal{BP}_{a,b,k,r}(m).\]
 Obviously, $\mathcal{BP}_{a,b,k,r}$ is the basis of $\mathcal{P}_{a,b,k,r}$. So, we have the following result. 
 \begin{thm}
$\mathcal{P}_{a,b,k,r}$ is a separable integer partition class with modulus $k$.
\end{thm}

We find that \eqref{gen-p-abkr-0} is equivalent to showing the following theorem.
\begin{thm}\label{add-thm-new*}
For $m\geq 1$,
\begin{equation*}
\begin{split}
&\quad\sum_{\pi\in\mathcal{BP}_{a,b,k,r}(m)}\mu^{\ell_{k,a}(\pi)}\nu^{\ell_{k,b}(\pi)}q^{|\pi|}\\
&=\mu^mq^{ma}+\sum_{s\geq 1}\sum_{h\geq 0}\mu^{m-h-s}\nu^{h+s}q^{(m-h-s)a+(h+s)b+k(s^2-s)}g_{k,r}(h,s){{m-h-s+1}\brack {s}}_{k}.
\end{split}
\end{equation*}
\end{thm}

In order to prove Theorem \ref{add-thm-new*}, we need the following lemma.

\begin{lem}\label{gen-bp-0-lem}
For $m\geq1$, let  $\mathcal{BP}_{a,b,k,r}(m,a)$ (resp. $\mathcal{BP}_{a,b,k,r}(m,b)$) be the set of partitions in $\mathcal{BP}_{a,b,k,r}(m)$ with the smallest part being $a$ (resp. $b$). Then, we have
\begin{equation}\label{gen-bp-0-lem-a}
\begin{split}
&\quad\sum_{\pi\in\mathcal{BP}_{a,b,k,r}(m,a)}\mu^{\ell_{k,a}(\pi)}\nu^{\ell_{k,b}(\pi)}q^{|\pi|}\\
&=\mu^mq^{ma}+\sum_{s\geq 1}\sum_{h\geq 0}\mu^{m-h-s}\nu^{h+s}q^{(m-h-s)a+(h+s)b+k(s^2-s)}g_{k,r}(h,s){{m-h-s}\brack {s}}_{k},
\end{split}
\end{equation}
and
\begin{equation}\label{gen-bp-0-lem-b}
\begin{split}
&\quad\sum_{\pi\in\mathcal{BP}_{a,b,k,r}(m,b)}\mu^{\ell_{k,a}(\pi)}\nu^{\ell_{k,b}(\pi)}q^{|\pi|}\\
&=\sum_{s\geq 1}\sum_{h\geq 0}\mu^{m-h-s}\nu^{h+s}q^{(m-h-s)a+(h+s)b+k(s^2-s)+k(m-h-2s+1)}g_{k,r}(h,s){{m-h-s}\brack {s-1}}_{k}.
\end{split}
\end{equation}
\end{lem}

\pf We first show \eqref{gen-bp-0-lem-a}. It is easy to see that there is only one partition
\[(\underbrace{a,a,\ldots,a}_{m'\text{s}})\]
in $\mathcal{BP}_{a,b,k,r}(m,a)$ without parts congruent to $b$ modulo $k$. For $s\geq 1$ and $h\geq 0$, let $\mathcal{BP}_{a,b,k,r}(m,a,h,s)$ be the set of partitions $\pi$ in $\mathcal{BP}_{a,b,k,r}(m,a)$ with $L_{k,b}(\pi)=k(s-1)+b$ and $h=\ell_{k,b}(\pi)-s$, where $L_{k,b}(\pi)$ is the largest part congruent to $b$ modulo $k$ in $\pi$. It remains to show that
 \begin{equation}\label{gen-bp-0-lem-a-proof-1}
 \begin{split}
&\quad\sum_{\pi\in\mathcal{BP}_{a,b,k,r}(m,a,h,s)}\mu^{\ell_{k,a}(\pi)}\nu^{\ell_{k,b}(\pi)}q^{|\pi|}\\
&=\mu^{m-h-s}\nu^{h+s}q^{(m-h-s)a+(h+s)b+k(s^2-s)}g_{k,r}(h,s){{m-h-s}\brack {s}}_{k}.
\end{split}
\end{equation}

Let $\pi$ be a partition in $\mathcal{BP}_{a,b,k,r}(m,a,h,s)$. It follows from $L_{k,b}(\pi)=k(s-1)+b$ that there exist parts
\[a,b,k+a,k+b,\ldots,k(s-1)+a,k(s-1)+b\]
in $\pi$. We remove one $a$, one $b$, one $k+a$, one $k+b$, \ldots, one $k(s-1)+a$ and one $k(s-1)+b$ from $\pi$ and denote the resulting partition by $\lambda$. Set $\lambda^{a}$ (resp. $\lambda^{b}$) to be the partition consisting of the parts congruent to $a$ (resp. $b$) modulo $k$ in $\lambda$. Under the condition that $h=\ell_{k,b}(\pi)-s$, we know that $\lambda^{b}$ is a partition in $\mathcal{G}_{b,k,r}(h,s)$, and so there are $m-2s-h$ parts  in $\lambda^{a}$ which do not exceed $ks+a$. By Theorem \ref{useful-gen-h-s}, we get
\begin{align*}
&\quad\sum_{\pi\in\mathcal{BP}_{a,b,k,r}(m,a,h,s)}\mu^{\ell_{k,a}(\pi)}\nu^{\ell_{k,b}(\pi)}q^{|\pi|}\\
&=\mu^sq^{k(s^2-s)/2+sa}\cdot\nu^sq^{k(s^2-s)/2+sb}\cdot\nu^hq^{hb}g_{k,r}(h,s)\cdot\mu^{m-2s-h}q^{(m-2s-h)a}{{m-h-s}\brack {s}}_{k}.
\end{align*}
We arrive at \eqref{gen-bp-0-lem-a-proof-1}, and thus \eqref{gen-bp-0-lem-a} holds. With a similar argument above, we can show that \eqref{gen-bp-0-lem-b} is valid.  The proof is complete.   \qed

Now, we are in a position to show Theorem \ref{add-thm-new*}.

{\noindent \bf Proof of Theorem \ref{add-thm-new*}.} In light of Lemma \ref{gen-bp-0-lem}, we have
\begin{align*}
&\quad\sum_{\pi\in\mathcal{BP}_{a,b,k,r}(m)}\mu^{\ell_{k,a}(\pi)}\nu^{\ell_{k,b}(\pi)}q^{|\pi|}\\
&=\sum_{\pi\in\mathcal{BP}_{a,b,k,r}(m,a)}\mu^{\ell_{k,a}(\pi)}\nu^{\ell_{k,b}(\pi)}q^{|\pi|}+\sum_{\pi\in\mathcal{BP}_{a,b,k,r}(m,b)}\mu^{\ell_{k,a}(\pi)}\nu^{\ell_{k,b}(\pi)}q^{|\pi|}\\
&=\mu^mq^{ma}+\sum_{s\geq 1}\sum_{h\geq 0}\mu^{m-h-s}\nu^{h+s}q^{(m-h-s)a+(h+s)b+k(s^2-s)}g_{k,r}(h,s)\\
&\qquad\qquad\qquad\qquad\times\left({{m-h-s}\brack {s}}_{k}+q^{k(m-h-2s+1)}{{m-h-s}\brack {s-1}}_{k}\right)\\
&=\mu^mq^{ma}+\sum_{s\geq 1}\sum_{h\geq 0}\mu^{m-h-s}\nu^{h+s}q^{(m-h-s)a+(h+s)b+k(s^2-s)}g_{k,r}(h,s){{m-h-s+1}\brack {s}}_{k},
\end{align*}
where the final equation follows from the following recurrence for the $q$-binomial coefficients \cite[(3.3.3)]{Andrews-1976} with $A=m-h-s+1$, $B=s$ and $q\rightarrow q^k$.
\begin{equation*}\label{bin-new-r-1}
{A\brack B}_1={{A-1}\brack{B}}_1+q^{A-B}{{A-1}\brack{B-1}}_1.
\end{equation*}
This competes the proof. \qed

\subsection{Proof of \eqref{gen-pp-abkr-0}}

The proof of \eqref{gen-pp-abkr-0} is similar to that of \eqref{gen-p-abkr-0}.
For $m\geq 1$, let $\mathcal{BP}'_{a,b,k,r}(m)$ denote the set of partitions $\pi=(\pi_1,\pi_2,\ldots,\pi_m)$ satisfying the conditions (1)-(3) in Definition \ref{basis-p-m}, and  for $1\leq i\leq m-r$, if $\pi_{i+1}=\pi_{i+2}=\cdots=\pi_{i+r}\equiv a\pmod{k}$ then $\pi_i\equiv b\pmod{k}$.

For example, there are thirteen partitions in $\mathcal{BP}'_{a,b,k,2}(4)$.
\[(b,b,a,a),(k+a,b,a,a),(b,b,b,a),(k+a,b,b,a),\]
\[(k+a,k+a,b,a),(k+b,k+a,b,a),(b,b,b,b),(k+a,b,b,b),(k+a,k+a,b,b),\]
\[(k+b,k+a,b,b),(k+b,k+a,k+a,b),(k+b,k+b,k+a,b),(2k+a,k+b,k+a,b).\]

Set
\[\mathcal{BP}'_{a,b,k,r}=\bigcup_{m\geq1}\mathcal{BP}'_{a,b,k,r}(m).\]
 Obviously, $\mathcal{BP}'_{a,b,k,r}$ is the basis of $\mathcal{P}'_{a,b,k,r}$. So, we have the following result.
 \begin{thm}
$\mathcal{P}'_{a,b,k,r}$ is a separable integer partition class with modulus $k$.
\end{thm}

We find that in order to prove \eqref{gen-pp-abkr-0}, it suffices to show the following lemma.

\begin{lem}\label{gen-bpp-0-lem}
For $m\geq1$, let  $\mathcal{BP}'_{a,b,k,r}(m,a)$ (resp. $\mathcal{BP}'_{a,b,k,r}(m,b)$) be the set of partitions in $\mathcal{BP}'_{a,b,k,r}(m)$ with the smallest part being $a$ (resp. $b$).  Then, we have
\begin{equation}\label{gen-bpp-0-lem-a}
\begin{split}
&\quad\sum_{\pi\in\mathcal{BP}'_{a,b,k,r}(m,a)}\mu^{\ell_{k,a}(\pi)}\nu^{\ell_{k,b}(\pi)}q^{|\pi|}\\
&=\sum_{s\geq 1}\sum_{h\geq 0}\mu^{h+s}\nu^{m-h-s}q^{(h+s)a+(m-h-s)b+k(s-1)^2}g_{k,r}(h,s){{m-h-s}\brack {s-1}}_{k},
\end{split}
\end{equation}
and
\begin{equation}\label{gen-bpp-0-lem-b}
\begin{split}
&\quad\sum_{\pi\in\mathcal{BP}'_{a,b,k,r}(m,b)}\mu^{\ell_{k,a}(\pi)}\nu^{\ell_{k,b}(\pi)}q^{|\pi|}\\
&=\nu^mq^{mb}+\sum_{s\geq 1}\sum_{h\geq 0}\mu^{h+s}\nu^{m-h-s}q^{(h+s)a+(m-h-s)b+ks^2+kh}g_{k,r}(h,s){{m-h-s}\brack {s}}_{k}.
\end{split}
\end{equation}
\end{lem}

\pf For a partition $\pi$ in $\mathcal{BP}_{a,b,k,r}(m,b)$, if we subtract $k-b+a$ from each of the part congruent to $a$ modulo $k$ in $\pi$ and subtract $b-a$ from each of the part congruent to $b$ modulo $k$ in $\pi$, then we can get a partition in $\mathcal{BP}'_{a,b,k,r}(m,a)$, and vice versa. Using \eqref{gen-bp-0-lem-b}, we have
\begin{align*}
\sum_{\pi\in\mathcal{BP}'_{a,b,k,r}(m,a)}\mu^{\ell_{k,a}(\pi)}\nu^{\ell_{k,b}(\pi)}q^{|\pi|}&=\sum_{s\geq 1}\sum_{h\geq 0}\nu^{m-h-s}\mu^{h+s}q^{-(m-h-s)(k-b+a)-(h+s)(b-a)}\\
&\quad\times q^{(m-h-s)a+(h+s)b+k(s^2-s)+k(m-h-2s+1)}g_{k,r}(h,s){{m-h-s}\brack {s-1}}_{k}.
\end{align*}
So, \eqref{gen-bpp-0-lem-a} is satisfied.

For a partition $\pi$ in $\mathcal{BP}_{a,b,k,r}(m,a)$, if we add $b-a$ to each of the part congruent to $a$ modulo $k$ in $\pi$ and add $k-b+a$ to each of the part congruent to $b$ modulo $k$ in $\pi$, then we can get a partition in $\mathcal{BP}'_{a,b,k,r}(m,b)$, and vice versa. Appealing to \eqref{gen-bp-0-lem-a},  we have
\begin{align*}
&\quad \sum_{\pi\in\mathcal{BP}'_{a,b,k,r}(m,b)}\mu^{\ell_{k,a}(\pi)}\nu^{\ell_{k,b}(\pi)}q^{|\pi|}=\nu^mq^{m(b-a)+ma}+\sum_{s\geq 1}\sum_{h\geq 0}\nu^{m-h-s}\mu^{h+s}\\
&\qquad\qquad\qquad\times q^{(m-h-s)(b-a)+(h+s)(k-b+a)}\cdot q^{(m-h-s)a+(h+s)b+k(s^2-s)}g_{k,r}(h,s){{m-h-s}\brack {s}}_{k}.
\end{align*}
So, \eqref{gen-bpp-0-lem-b} is valid.  The proof is complete.  \qed

\section{Proof of Theorem \ref{thm-r-abkr-0}}

In this section, we assume that $k,a,b,c$ and $r$ are integers such that  $k\geq c>b>a\geq 1$ and $r\geq 1$. The objective of this section is to give a proof of Theorem \ref{thm-r-abkr-0}. We first introduce the following two sets which are the main ingredients in the construction of the basis of $\mathcal{R}_{a,b,c,k}$ and $\mathcal{R}_{a,b,c,k,r}$.

\begin{defi}\label{basis-r-m}
For $m\geq 1$, let $\mathcal{BR}_{a,b,c,k}(m)$ (resp. $\mathcal{BR}_{a,b,c,k,r}(m)$) be the set of partitions $\pi=(\pi_1,\pi_2,\ldots,\pi_m)$ in  $\mathcal{R}_{a,b,c,k}$ (resp. $\mathcal{R}_{a,b,c,k,r}$) such that
 \begin{itemize}
 \item[(1)] $\pi_m=a$ or $b$ or $c$;
\item[(2)] for $1\leq i<m$, $\pi_i-\pi_{i+1}\leq k$ with strict inequality if $\pi_{i+1}$ is congruent to $a$ or $b$ modulo $k$.
\end{itemize}
\end{defi}

It is clear that for $m\geq 1$, the number of partitions in $\mathcal{BR}_{a,b,c,k}(m)$ does not exceed $3^m$.
By definition, we know that for $m\geq 1$, $\mathcal{BR}_{a,b,c,k,r}(m)$ is the set of partitions $\pi=(\pi_1,\pi_2,\ldots,\pi_m)$ in $\mathcal{BR}_{a,b,c,k}(m)$ such that for $1\leq i\leq m-r+1$, if $\pi_{i+1}=\pi_{i+2}=\cdots=\pi_{i+r-1}\equiv b\pmod{k}$ then $\pi_i\equiv c\pmod{k}$, and so $\mathcal{BR}_{a,b,c,k,r}(m)$ is a subset of $\mathcal{BR}_{a,b,c,k}(m)$.

For example, there are seventeen partitions in $\mathcal{BR}_{a,b,c,k,2}(3)$.
\[(a,a,a),(b,a,a),(c,a,a),(c,b,a),(k+a,c,a),(k+b,c,a),\]
\[(k+c,c,a),(k+a,c,b),(k+b,c,b),(k+c,c,b),(k+a,k+a,c),(k+b,k+a,c),\]
\[(k+c,k+a,c),(k+c,k+b,c),(2k+a,k+c,c),(2k+b,k+c,c),(2k+c,k+c,c).\]
Besides, there are four partitions belonging to $\mathcal{BR}_{a,b,c,k}(3)$ but not in  $\mathcal{BR}_{a,b,c,k,2}(3)$.
\[(b,b,a),(b,b,b),(c,b,b),(k+b,k+b,c).\]

Set
\[\mathcal{BR}_{a,b,c,k}=\bigcup_{m\geq1}\mathcal{BR}_{a,b,c,k}(m),\]
and
\[\mathcal{BR}_{a,b,c,k,r}=\bigcup_{m\geq1}\mathcal{BR}_{a,b,c,k,r}(m).\]
Obviously, $\mathcal{BR}_{a,b,c,k}$ and $\mathcal{BR}_{a,b,c,k,r}$ are the basis of $\mathcal{R}_{a,b,c,k}$ and $\mathcal{R}_{a,b,c,k,r}$ respectively. So, we have the following result. 
 \begin{thm}
$\mathcal{R}_{a,b,c,k}$ and $\mathcal{R}_{a,b,c,k,r}$ are separable integer partition classes with modulus $k$.
\end{thm}

We find that in order to prove Theorem \ref{thm-r-abkr-0}, it suffices to show the following lemma.
\begin{lem}For $m\geq 1$,
\begin{equation}\label{gen-basis-r-abkr-0}
 \begin{split}
&\quad\sum_{\pi\in\mathcal{BR}_{a,b,c,k}(m)}\mu^{\ell_{k,a}(\pi)}\nu^{\ell_{k,b}(\pi)}\omega^{\ell_{k,c}(\pi)}q^{|\pi|}\\
&=\sum_{s,h\geq 0}\mu^{m-h-s}\nu^{h}\omega^sq^{(m-h-s)a+hb+sc+k(s^2-s)/2}{{h+s}\brack {s}}_{k}{{m-h}\brack {s}}_{k},
 \end{split}
 \end{equation}
 and
\begin{equation}\label{gen-basis-rr-abkr-0}
\begin{split}
&\quad\sum_{\pi\in\mathcal{BR}_{a,b,c,k,r}(m)}\mu^{\ell_{k,a}(\pi)}\nu^{\ell_{k,b}(\pi)}\omega^{\ell_{k,c}(\pi)}q^{|\pi|}\\
&=\sum_{s,h\geq 0}\mu^{m-h-s}\nu^{h}\omega^sq^{(m-h-s)a+hb+sc+k(s^2-s)/2}g_{k,r}(h,s+1){{m-h}\brack {s}}_{k}.
\end{split}
\end{equation}
\end{lem}

\pf For $s\geq 1$ and $h\geq 0$, let $\mathcal{BR}_{a,b,c,k}(m,h,s)$ (resp. $\mathcal{BR}_{a,b,c,k,r}(m,h,s)$) be the set of partitions in $\mathcal{BR}_{a,b,c,k}(m)$ (resp. $\mathcal{BR}_{a,b,c,k,r}(m)$) such that there are $h$ parts congruent to $b$ modulo $k$ and the largest part congruent $c$ modulo $k$ is $k(s-1)+c$.

For  a partition $\pi$  in $\mathcal{BR}_{a,b,c,k}(m,h,s)$, let $\pi^{a}$, $\pi^{b}$ and $\pi^{c}$ be the partitions consisting of the parts congruent to $a$, $b$ and $c$ modulo $k$ in $\pi$ respectively. By definition, we see that
\begin{itemize}
\item[(1)] $\pi^c=(k(s-1)+c,\ldots,k+c,c)$;
\item[(2)] there are $h$ parts in $\pi^{b}$ which do not exceed $ks+b$;
\item[(3)] there are $m-h-s$ parts  in $\pi^{a}$ which do not exceed $ks+a$.
\end{itemize}
 Moreover, $\pi$ is also a partition in $\mathcal{BR}_{a,b,c,k,r}(m,h,s)$ if and only if each part of  $\pi^b$ appears at most $r-1$ times. Namely, $\pi$ is a partition in $\mathcal{BR}_{a,b,c,k,r}(m,h,s)$ if and only if $\pi^b\in\mathcal{G}_{b,k,r}(h,s+1)$. So, we have
 \begin{equation}\label{gen-basis-r-abkrhs-0}
 \begin{split}
&\quad\sum_{\pi\in\mathcal{BR}_{a,b,c,k}(m,h,s)}\mu^{\ell_{k,a}(\pi)}\nu^{\ell_{k,b}(\pi)}\omega^{\ell_{k,c}(\pi)}q^{|\pi|}\\
&=\omega^sq^{k(s^2-s)/2+sc}\cdot\nu^{h}q^{hb}{{h+s}\brack {s}}_{k}\cdot\mu^{m-h-s}q^{(m-h-s)a}{{m-h}\brack {s}}_{k}\\
&=\mu^{m-h-s}\nu^{h}\omega^sq^{(m-h-s)a+hb+sc+k(s^2-s)/2}{{h+s}\brack {s}}_{k}{{m-h}\brack {s}}_{k},
 \end{split}
 \end{equation}
 and by Theorem \ref{useful-gen-h-s}, we get
\begin{equation}\label{gen-basis-rr-abkrhs-0}
\begin{split}
&\quad\sum_{\pi\in\mathcal{BR}_{a,b,c,k,r}(m,h,s)}\mu^{\ell_{k,a}(\pi)}\nu^{\ell_{k,b}(\pi)}\omega^{\ell_{k,c}(\pi)}q^{|\pi|}\\
&=\omega^sq^{k(s^2-s)/2+sc}\cdot \nu^hq^{hb}g_{k,r}(h,s+1)\cdot\mu^{m-h-s}q^{(m-h-s)a}{{m-h}\brack {s}}_{k}\\
&=\mu^{m-h-s}\nu^{h}\omega^sq^{(m-h-s)a+hb+sc+k(s^2-s)/2}g_{k,r}(h,s+1){{m-h}\brack {s}}_{k}.
\end{split}
\end{equation}

It is easy to see that the set of partitions in $\mathcal{BR}_{a,b,c,k}(m)$ and $\mathcal{BR}_{a,b,c,k,r}(m)$ without parts congruent to $c$ modulo $k$ are
\[\{(\underbrace{b,\ldots,b}_{h'\text{s}},\underbrace{a,\ldots,a}_{(m-h)'\text{s}})\mid 0\leq h\leq m\},\]
and
\[\{(\underbrace{b,\ldots,b}_{h'\text{s}},\underbrace{a,\ldots,a}_{(m-h)'\text{s}})\mid 0\leq h\leq \text{min}\{m,r-1\}\}.\]

So, we get
\begin{align*}
&\quad\sum_{\pi\in\mathcal{BR}_{a,b,c,k}(m)}\mu^{\ell_{k,a}(\pi)}\nu^{\ell_{k,b}(\pi)}\omega^{\ell_{k,c}(\pi)}q^{|\pi|}\\
&=\sum_{h=0}^{m}\mu^{m-h}\nu^{h}q^{(m-h)a+hb}+\sum_{s\geq 1}\sum_{h\geq 0}\sum_{\pi\in\mathcal{BR}_{a,b,c,k}(m,h,s)}\mu^{\ell_{a}(\pi)}\nu^{\ell_{b}(\pi)}\omega^{\ell_{c}(\pi)}q^{|\pi|},
\end{align*}
and
\begin{align*}
&\quad\sum_{\pi\in\mathcal{BR}_{a,b,c,k,r}(m)}\mu^{\ell_{k,a}(\pi)}\nu^{\ell_{k,b}(\pi)}\omega^{\ell_{k,c}(\pi)}q^{|\pi|}\\
&=\sum_{h=0}^{\text{min}\{m,r-1\}}\mu^{m-h}\nu^{h}q^{(m-h)a+hb}+\sum_{s\geq 1}\sum_{h\geq 0}\sum_{\pi\in\mathcal{BR}_{a,b,c,k,r}(m,h,s)}\mu^{\ell_{a}(\pi)}\nu^{\ell_{b}(\pi)}\omega^{\ell_{c}(\pi)}q^{|\pi|}.
\end{align*}
Combining with \eqref{gen-basis-r-abkrhs-0} and \eqref{gen-basis-rr-abkrhs-0}, we arrive at  \eqref{gen-basis-r-abkr-0} and \eqref{gen-basis-rr-abkr-0}. This completes the proof.  \qed

\section{Proof of Theorem \ref{main-over}}

In this section, we aim to give a proof of Theorem \ref{main-over}. To do this, we will show \eqref{eqn-over-f-over} and \eqref{eqn-over-f-non-over} in Section 5.1 and show \eqref{eqn-over-L-over} and \eqref{eqn-over-L-non-over} in Section 5.2.

\subsection{Proofs of \eqref{eqn-over-f-over} and \eqref{eqn-over-f-non-over}}

In this subsection, an overpartition is defined as a partition in which the first occurrence of a number may be overlined. We first give a proof of \eqref{eqn-over-f-over}. For $m\geq 1$, let  $\overline{\mathcal{BF}}(m)$ be the set of overpartitions $\pi=(\pi_1,\pi_2,\ldots,\pi_m)$ such that
\begin{itemize}
\item[(1)] $\pi_m=1$ or $\overline{1}$;
\item[(2)] for $1\leq i<m$,
\begin{itemize}
\item[(2.1)] if $\pi_{i+1}$ is non-overlined, then $|\pi_i|=|\pi_{i+1}|$;
\item[(2.2)] if $\pi_{i+1}$ is overlined, then $|\pi_i|=|\pi_{i+1}|+1$ and $\pi_i$ is non-overlined.
\end{itemize}
\end{itemize}

Clearly, the number of overpartitions  in $\overline{\mathcal{BF}}(m)$ does not exceed $2^m$. For example, there are eight overpartitions in $\overline{\mathcal{BF}}(4)$.
\[(1,1,1,1),(\overline{1},1,1,1),(2,\overline{1},1,1),(2,2,\overline{1},1),(\overline{2},2,\overline{1},1),(2,2,2,\overline{1}),(\overline{2},2,2,\overline{1}),(3,\overline{2},2,\overline{1}).\]

Set
\[\overline{\mathcal{BF}}=\bigcup_{m\geq 1}\overline{\mathcal{BF}}(m).\]
Obviously, $\overline{\mathcal{BF}}$ is the basis of $\overline{\mathcal{F}}$. So, we have the following result.
 \begin{thm}
$\overline{\mathcal{F}}$ is a separable overpartition class.
\end{thm}

Now, we are in a position to show \eqref{eqn-over-f-over}.

{\noindent\bf Proof of \eqref{eqn-over-f-over}.} For $m\geq 1$ and $s\geq 0$, let $\overline{\mathcal{BF}}(m,s)$ be the set of overpartitions in $\overline{\mathcal{BF}}(m)$ with $s$ overlined parts. Then, it suffices to show that
\begin{equation}\label{proof-F-over-ms}
\sum_{\pi\in\overline{\mathcal{BF}}(m,s)}q^{|\pi|}=q^{m+s^2-s}{{m-s+1}\brack{s}}_1.
\end{equation}
We consider the following two cases.

Case 1: For $s=0$, it is easy to see that there is only one overpartition
\[(\underbrace{1,1,\ldots,1}_{m'\text{s}})\]
in $\overline{\mathcal{BF}}(m,0)$. This implies that
\[\sum_{\pi\in\overline{\mathcal{BF}}(m,0)}q^{|\pi|}=q^{m},\]
which agrees with \eqref{proof-F-over-ms} for $s=0$.

Case 2: For $s\geq 1$, let $\pi$ be an overpartition in $\overline{\mathcal{BF}}(m,s)$. Then, there exist parts
\[\overline{1},2,\overline{2},\ldots,s,\overline{s}\]
in $\pi$. We remove one $\overline{1}$, one $2$, one $\overline{2}$, \ldots, one $s$ and one $\overline{s}$ from $\pi$ and denote the resulting overpartition by $\lambda$. It is clear that there are no overlined parts in $\lambda$ and there are $m-2s+1$ non-overlined parts in $\lambda$ which do not exceed $s+1$. So, we get
\[\sum_{\pi\in\overline{\mathcal{BF}}(m,s)}q^{|\pi|}=q^{s^2+s-1}\cdot q^{m-2s+1}{{m-s+1}\brack{s}}_1=q^{m+s^2-s}{{m-s+1}\brack{s}}_1.\]
We arrive at \eqref{proof-F-over-ms} for $s\geq 1$. This completes the proof.  \qed

Then, we proceed to show \eqref{eqn-over-f-non-over}. In the remaining of this subsection, we fix $r\geq 1$. For $m\geq 1$, let  ${\mathcal{BF}}_r(m)$ be the set of overpartitions $\pi=(\pi_1,\pi_2,\ldots,\pi_m)$ such that
\begin{itemize}
\item[(1)] $\pi_m=1$ or $\overline{1}$;
\item[(2)] for $1\leq i<m$,
\begin{itemize}
\item[(2.1)] if $\pi_{i+1}$ is non-overlined, then $|\pi_i|=|\pi_{i+1}|$;
\item[(2.2)] if $\pi_{i+1}$ is overlined, then $|\pi_i|=|\pi_{i+1}|+1$;
\end{itemize}
\item[(3)] for $1\leq i\leq m-r+1$, if $\pi_{i+1},\pi_{i+2},\ldots,\pi_{i+r-1}$ are non-overlined then $\pi_i$ is overlined.
\end{itemize}

Clearly, the number of overpartitions  in ${\mathcal{BF}}_r(m)$ does not exceed $2^m$. For example, there are thirteen overpartitions in ${\mathcal{BF}}_3(4)$.
\[(2,\overline{1},1,1),(\overline{2},\overline{1},1,1),(2,2,\overline{1},1),(\overline{2},2,\overline{1},1),(3,\overline{2},\overline{1},1),(\overline{3},\overline{2},\overline{1},1),\]
\[(\overline{2},2,2,\overline{1}),(3,\overline{2},2,\overline{1}),(\overline{3},\overline{2},2,\overline{1}),(3,3,\overline{2},\overline{1}),(\overline{3},3,\overline{2},\overline{1}),(4,\overline{3},\overline{2},\overline{1}),(\overline{4},\overline{3},\overline{2},\overline{1}).\]

Set
\[{\mathcal{BF}}_r=\bigcup_{m\geq 1}{\mathcal{BF}}_r(m).\]
Obviously, ${\mathcal{BF}}_r$ is the basis of ${\mathcal{F}}_r$. So, we have the following result.
 \begin{thm}
${\mathcal{F}}_r$ is a separable overpartition class.
\end{thm}

We conclude this subsection with a proof of \eqref{eqn-over-f-non-over}.

{\noindent\bf Proof of \eqref{eqn-over-f-non-over}.}  For $m\geq 1$ and $s\geq 0$, let  ${\mathcal{BF}}_r(m,s)$ be the set of overpartitions in ${\mathcal{BF}}_r(m)$ with $s$ overlined parts. Then, it suffices to show that
\begin{equation}\label{proof-F-non-over-rms}
\sum_{\pi\in{\mathcal{BF}}_r(m,s)}q^{|\pi|}=q^{m+(s^2-s)/2}g_{1,r}(m-s,s+1).
\end{equation}
We consider the following two cases.

Case 1: For $s=0$, the right-hand side of \eqref{proof-F-non-over-rms} is $q^{m}g_{1,r}(m,1)$. It follows from Theorem \ref{useful-gen-h-s} that $q^{m}g_{1,r}(m,1)$ is the generating function for the partitions in
\[\mathcal{G}_{1,1,r}(m,1)=\left\{\begin{array}{ll}\{(\underbrace{1,1,\ldots,1}_{m'\text{s}})\},&\text{if }0\leq m\leq r-1,\\
\emptyset,&\text{if }m>r-1.
\end{array}\right.\]
Clearly, we have $\mathcal{G}_{1,1,r}(m,1)=\mathcal{BF}_r(m,0)$. This implies that  \eqref{proof-F-non-over-rms} holds for $s=0$.

Case 2: For $s\geq 1$,  let $\pi$ be an overpartition in ${\mathcal{BF}}_r(m,s)$. Then, there exist parts
\[\overline{1},\overline{2},\ldots,\overline{s}\] in $\pi$. We remove one $\overline{1}$, one $\overline{2}$, \ldots, one $\overline{s}$ from $\pi$ and denote the resulting overpartition by $\lambda$. Clearly, there are no overlined parts in $\lambda$, there are $m-s$ non-overlined parts in $\lambda$ which do not exceed $s+1$, and each part of $\lambda$ appears at most $r-1$ times. It yields that $\lambda$ is a partition in $\mathcal{G}_{1,1,r}(m-s,s+1)$. Utilizing Theorem \ref{useful-gen-h-s}, we get
\[\sum_{\pi\in\mathcal{BF}_r(m,s)}q^{|\pi|}=q^{(s^2+s)/2}\cdot q^{m-s}g_{1,r}(m-s,s+1)=q^{m+(s^2-s)/2}g_{1,r}(m-s,s+1).\]
We arrive at \eqref{proof-F-non-over-rms} for $s\geq 1$. This completes the proof.  \qed

\subsection{Proofs of \eqref{eqn-over-L-over} and \eqref{eqn-over-L-non-over}}

In this subsection, an overpartition is defined as a partition in which the last occurrence of a number may be overlined. We first give a proof of \eqref{eqn-over-L-over}. For $m\geq 1$, let  $\overline{\mathcal{BL}}(m)$ be the set of overpartitions $\pi=(\pi_1,\pi_2,\ldots,\pi_m)$ such that
\begin{itemize}
\item[(1)] $\pi_m=1$ or $\overline{1}$;
\item[(2)] for $1\leq i<m$,
\begin{itemize}
\item[(2.1)] if $\pi_{i}$ is non-overlined, then $|\pi_i|=|\pi_{i+1}|$;
\item[(2.2)] if $\pi_{i}$ is overlined, then $|\pi_i|=|\pi_{i+1}|+1$ and $\pi_{i+1}$ is non-overlined.
\end{itemize}
\end{itemize}

Clearly, the number of overpartitions  in $\overline{\mathcal{BL}}(m)$ does not exceed $2^m$. For example, there are eight overpartitions in $\overline{\mathcal{BL}}(4)$.
\[(1,1,1,1),(\overline{2},1,1,1),(2,\overline{2},1,1),(2,2,\overline{2},1),(\overline{3},2,\overline{2},1),(1,1,1,\overline{1}),(\overline{2},1,1,\overline{1}),(2,\overline{2},1,\overline{1}).\]

Set
\[\overline{\mathcal{BL}}=\bigcup_{m\geq 1}\overline{\mathcal{BL}}(m).\]
Obviously, $\overline{\mathcal{BL}}$ is the basis of $\overline{\mathcal{L}}$. So, we have the following result.
 \begin{thm}
$\overline{\mathcal{L}}$ is a separable overpartition class.
\end{thm}

Now, we are in a position to show \eqref{eqn-over-L-over}.

{\noindent\bf Proof of \eqref{eqn-over-L-over}.} For $m\geq 1$ and $s\geq 0$, let $\overline{\mathcal{BL}}(m,s)$ be the set of overpartitions in $\overline{\mathcal{BL}}(m)$ with $s$ overlined parts. Then, it suffices to show that
\begin{equation}\label{proof-L-over-ms}
\sum_{\pi\in\overline{\mathcal{BL}}(m,s)}q^{|\pi|}=q^{m+(s-1)^2}\left({{m-s}\brack{s-1}}_1+q^{2s-1}{{m-s}\brack{s}}_1\right).
\end{equation}
We consider the following two cases.

Case 1: For $s=0$, the right-hand side of \eqref{proof-L-over-ms} equals $q^m$, which is the generating function for the overpartitions in
\[\overline{\mathcal{BL}}(m,0)=\{(\underbrace{1,1,\ldots,1}_{m'\text{s}})\}.\]
So, \eqref{proof-L-over-ms} holds for $s=0$.

Case 2: For $s\geq 1$, let $\overline{\mathcal{BL}}'(m,s)$ (resp. $\overline{\mathcal{BL}}''(m,s)$) be the set of overpartitions in $\overline{\mathcal{BL}}(m,s)$ with the smallest part being $\overline{1}$ (resp. $1$).  For an overpartition  $\pi$ in $\overline{\mathcal{BL}}'(m,s)$, it is easy to see that  there exist parts
\[\overline{1},1,,\ldots,\overline{s-1},s-1,\overline{s}\]
in $\pi$.  We remove one $\overline{1}$, one $1$, \ldots, one $\overline{s-1}$, one $s-1$ and one $\overline{s}$ from $\pi$ and denote the resulting overpartition by $\lambda$. It is clear that there are no overlined parts in $\lambda$ and there are $m-2s+1$ non-overlined parts in $\lambda$ which do not exceed $s$. So, we have
\begin{equation}\label{proof-BL'}
\sum_{\pi\in\overline{\mathcal{BL}}'(m,s)}q^{|\pi|}=q^{s^2}\cdot q^{m-2s+1}{{m-s}\brack{s-1}}_1=q^{m+(s-1)^2}{{m-s}\brack{s-1}}_1.
\end{equation}

For an overpartition $\alpha=(\alpha_1,\alpha_2,\ldots,\alpha_m,\alpha_{m+1})$ in $\overline{\mathcal{BL}}'(m+1,s+1)$, we have $\alpha_{m+1}=\overline{1}$ and $\alpha_{m}=1$. If we remove $\alpha_{m+1}=\overline{1}$ from $\alpha$, then we can get an overpartition in $\overline{\mathcal{BL}}''(m,s)$, and vice versa. This implies that
\begin{equation}\label{proof-BL''}
\sum_{\pi\in\overline{\mathcal{BL}}''(m,s)}q^{|\pi|}=q^{-1}\sum_{\pi\in\overline{\mathcal{BL}}'(m+1,s+1)}q^{|\pi|}=q^{m+s^2}{{m-s}\brack{s}}_1,
\end{equation}
where the final equation follows from \eqref{proof-BL'}. Combining \eqref{proof-BL'} and \eqref{proof-BL''}, we get
\[
\sum_{\pi\in\overline{\mathcal{BL}}(m,s)}q^{|\pi|}=\sum_{\pi\in\overline{\mathcal{BL}}'(m,s)}q^{|\pi|}+\sum_{\pi\in\overline{\mathcal{BL}}''(m,s)}q^{|\pi|}=q^{m+(s-1)^2}{{m-s}\brack{s-1}}_1+q^{m+s^2}{{m-s}\brack{s}}_1.
\]
We arrive at \eqref{proof-L-over-ms} for $s\geq 1$. The proof is complete.  \qed

Then, we turn to proving \eqref{eqn-over-L-non-over}. In the remaining of this subsection, we fix $r\geq 1$. For $m\geq 1$, let  ${\mathcal{BL}}_r(m)$ be the set of overpartitions $\pi=(\pi_1,\pi_2,\ldots,\pi_m)$ such that
\begin{itemize}
\item[(1)] $\pi_m=1$ or $\overline{1}$;
\item[(2)] for $1\leq i<m$,
\begin{itemize}
\item[(2.1)] if $\pi_{i}$ is non-overlined, then $|\pi_i|=|\pi_{i+1}|$;
\item[(2.2)] if $\pi_{i}$ is overlined, then $|\pi_i|=|\pi_{i+1}|+1$;
\end{itemize}
\item[(3)] for $1\leq i\leq m-r+1$, if $\pi_{i+1},\pi_{i+2},\ldots,\pi_{i+r-1}$ are non-overlined then $\pi_i$ is overlined.
\end{itemize}

Clearly, the number of overpartitions  in ${\mathcal{BL}}_r(m)$ does not exceed $2^m$. For example, there are thirteen overpartitions in ${\mathcal{BL}}_3(4)$.
\[(2,\overline{2},1,1),(\overline{3},\overline{2},1,1),(2,2,\overline{2},1),(\overline{3},2,\overline{2},1),(3,\overline{3},\overline{2},1),(\overline{4},\overline{3},\overline{2},1),\]
\[(\overline{2},1,1,\overline{1}),(2,\overline{2},1,\overline{1}),(\overline{3},\overline{2},1,\overline{1}),(2,2,\overline{2},\overline{1}),(\overline{3},2,\overline{2},\overline{1}),(3,\overline{3},\overline{2},\overline{1}),(\overline{4},\overline{3},\overline{2},\overline{1}).\]

Set
\[{\mathcal{BL}}_r=\bigcup_{m\geq 1}{\mathcal{BL}}_r(m).\]
Obviously, ${\mathcal{BL}}_r$ is the basis of ${\mathcal{L}}_r$. So, we have the following result. 
 \begin{thm}
${\mathcal{L}}_r$ is a separable overpartition class.
\end{thm}

Finally, we present a proof of \eqref{eqn-over-L-non-over}.

{\noindent\bf Proof of \eqref{eqn-over-L-non-over}.}  For $m\geq 1$ and $s\geq 0$, let  ${\mathcal{BL}}_r(m,s)$ be the set of overpartitions in ${\mathcal{BL}}_r(m)$ with $s$ overlined parts. By definition, we have
\[\mathcal{BL}_r(m,0)=\left\{\begin{array}{ll}\{(\underbrace{1,1,\ldots,1}_{m'\text{s}})\},&\text{if }0\leq m\leq r-1,\\
\emptyset,&\text{if }m>r-1,
\end{array}\right.\]
and so
\begin{equation*}\label{proof-bl-m-s-0}
\sum_{\pi\in\mathcal{BL}_r(m,0)}q^{|\pi|}=q^m{{r-m-1}\brack {0}}_{1}.
\end{equation*}

 It remains to show that for $s\geq 1$,
\begin{equation}\label{proof-L-non-over-rms}
\sum_{\pi\in{\mathcal{BL}}_r(m,s)}q^{|\pi|}=q^{m+(s^2-s)/2}\left(g_{1,r}(m-s,s)+q^m\sum_{j=1}^{r-1}q^{-j}g_{1,r}(m-j-s,s)\right).
\end{equation}

Let ${\mathcal{BL}}'_r(m,s)$ (resp. ${\mathcal{BL}}''_r(m,s)$) be the set of overpartitions in ${\mathcal{BL}}_r(m,s)$ with the smallest part being $\overline{1}$ (resp. $1$).  For an overpartition  $\pi$ in ${\mathcal{BL}}'_r(m,s)$, it is easy to see that  there exist parts
\[\overline{1},\overline{2},\ldots,\overline{s}\]
 in $\pi$. We remove one $\overline{1}$, one $\overline{2}$, \ldots, one $\overline{s}$ from $\pi$ and denote the resulting overpartition by $\lambda$.  Clearly, there are no overlined parts in $\lambda$, there are $m-s$ non-overlined parts in $\lambda$ which do not exceed $s$, and each part of $\lambda$ appears at most $r-1$ times. It yields that $\lambda$ is a partition in $\mathcal{G}_{1,1,r}(m-s,s)$.  Using Theorem \ref{useful-gen-h-s}, we get
\begin{equation}\label{final-use}
\sum_{\pi\in\mathcal{BL}'_r(m,s)}q^{|\pi|}=q^{(s^2+s)/2}\cdot q^{m-s}g_{1,r}(m-s,s)=q^{m+(s^2-s)/2}g_{1,r}(m-s,s).
\end{equation}

For $1\leq j\leq r-1$, assume that  $\alpha=(\alpha_1,\alpha_2,\ldots,\alpha_{m-j})$ is an overpartition in ${\mathcal{BL}}'_r(m-j,s)$, then \[(\alpha_1+1,\alpha_2+1,\ldots,\alpha_{m-j}+1,\underbrace{1,\ldots,1}_{j'\text{s}})\]
is an overpartition in ${\mathcal{BL}}''_r(m,s)$ with $j'$s $1$, and vice versa. This implies that

\begin{equation}\label{proof-BL''-r}
\sum_{\pi\in{\mathcal{BL}}''_r(m,s)}q^{|\pi|}=\sum_{j=1}^{r-1}q^{m}\sum_{\pi\in{\mathcal{BL}}'_r(m-j,s)}q^{|\pi|}=q^{2m+(s^2-s)/2}\sum_{j=1}^{r-1}q^{-j}g_{1,r}(m-j-s,s),
\end{equation}
where the final equation follows from \eqref{final-use}. Combining \eqref{final-use} and \eqref{proof-BL''-r},  we get
\begin{align*}
\sum_{\pi\in{\mathcal{BL}}_r(m,s)}q^{|\pi|}&=\sum_{\pi\in{\mathcal{BL}}'_r(m,s)}q^{|\pi|}+\sum_{\pi\in{\mathcal{BL}}''_r(m,s)}q^{|\pi|}\\
&=q^{m+(s^2-s)/2}g_{1,r}(m-s,s)+q^{2m+(s^2-s)/2}\sum_{j=1}^{r-1}q^{-j}g_{1,r}(m-j-s,s).
\end{align*}
We arrive at \eqref{proof-L-non-over-rms} for $s\geq 1$. The proof is complete.  \qed

\end{document}